\documentclass[12pt]{amsart}
\usepackage{amsmath}
\usepackage{amsfonts}
\usepackage{amsthm, upref}
\usepackage{graphicx}
\usepackage[usenames, dvipsnames]{color}

    \numberwithin{equation}{section}

\usepackage{bm}
         
         \bmdefine\alphab{\mathbf{\alpha}}
\bmdefine\betab{\mathbf{\beta}}
\bmdefine\sigmab{\mathbf{\sigma}}

\newcommand{\comment}[1]{}
\newcommand{\eq}{\begin{equation}}
\newcommand{\en}{\end{equation}}
\newcommand{\rr}{\mathbb{R}}
\newcommand{\nn}{\mathbb{N}}

\newcommand{\ep}{\hfill $\Box$}

\begin{document}

\theoremstyle{plain}
\newtheorem{thm}{Theorem}
\newtheorem{lemma}[thm]{Lemma}
\newtheorem{prop}[thm]{Proposition}
\newtheorem{cor}[thm]{Corollary}

\theoremstyle{definition}
\newtheorem{defn}{Definition}
\newtheorem{cond}{Condition}
\newtheorem{asmp}{Assumption}
\newtheorem{notn}{Notation}
\newtheorem{prb}{Problem}

\theoremstyle{remark}
\newtheorem{rmk}{Remark}
\newtheorem{exm}{Example}
\newtheorem{clm}{Claim}

\title[Strong solutions]{Strong solutions of stochastic 
equations with rank-based coefficients}

\author{Tomoyuki Ichiba, Ioannis Karatzas and Mykhaylo Shkolnikov}
\address{Department of Statistics and Applied Probability \\ University of California \\ Santa Barbara, CA 93106}
\email{ichiba@pstat.ucsb.edu}
\address{INTECH Investment Management \\ One Palmer Square\\ Princeton, NJ 08542 and Columbia University \\ Department of Mathematics\\ New York, NY 10027}
\email{ik@enhanced.com, ik@math.columbia.edu}
\address{INTECH Investment Management \\ One Palmer Square\\ Princeton, NJ 08542 and Stanford University \\ Department of Mathematics\\ Stanford, CA 94305}
\email{mshkolni@math.stanford.edu}

\keywords{Stochastic differential equations, strong existence, strong uniqueness, equations with rank-based coefficients, Brownian particles, triple collisions.}

    \subjclass[2000]{60H10, 60J60, 60J65}

\thanks{This research was partially supported by NSF grants DMS-08-06211 and DMS-09-05754.}

\date{\today}

\begin{abstract}
We study finite and countably infinite systems of stochastic differential equations, in which the drift and diffusion coefficients of each component  (particle) 
are determined by its rank in the vector of all components of the solution. We show that strong existence and uniqueness hold until the first time   three particles collide. Motivated by this result, we  improve significantly the existing conditions for the absence of such triple collisions in the case of finite-dimensional systems, and provide the first condition of this type   for  systems with a   countable infinity  of particles.     
\end{abstract}

\maketitle

\section{Introduction}

We study the following system of stochastic differential equations:
\eq\label{sdeI}
\mathrm{d}X_i(t)\,=\sum_{j\in I} \mathbf{1}_{\{X_i(t)=X_{(j)}(t)\}} \, \delta_j\, \mathrm{d}t 
+ \sum_{j\in I} \mathbf{1}_{\{X_i(t)=X_{(j)}(t)\}}\, \sigma_j\,\mathrm{d} W_i(t) 
\en
for $\, i \in I\,$. Here $I=\{1,\dots,n\}$ for some $\,n\in\nn\,$, or $I=\nn\,$; $\delta_j\,$, $j\in I$ are real constants; $\sigma_j$, $j\in I$ are strictly positive real constants; $(W_i:\;i\in I)$ is a system of independent standard Brownian motions; and 
\eq
X_{(1)}(t)\leq X_{(2)}(t)\leq X_{(3)}(t)\leq\ldots 
\en
is the ordered particle configuration at time $\,t\,$. In addition, we let the initial configuration be deterministic and satisfy  
\eq
X_1(0)<X_2(0)<X_3(0)<\ldots\,.
\en
Ties in the ordered particle configuration are resolved in accordance with the initial ranking of particles; for instance, we set $\, X_{(i)} (t) = X_i (t)\,$, $\, i=1, \ldots, n\,$ whenever $\, X_1 (t) = \ldots =  X_n (t)\,$.  
We shall write $X$ for $(X_i:\;i\in I)$ and $W$ for $(W_i:\;i\in I)$.

\medskip

In the case $I=\nn$ we work under the following assumption.

\asmp{\label{asmp} {\it If $I=\nn$, we assume that there is an integer $M\in\nn$ such that we have 
\begin{eqnarray}
&&\delta_M=\delta_{M+1}=\ldots\,\,,\\
&&\sigma_M=\sigma_{M+1}=\ldots\,\,.
\end{eqnarray}
Moreover, we assume that there exist constants $\gamma_1>0$, $\gamma_2\in\rr\,$ such that
\eq
\label{lin_growth}
X_i(0)\, \geq  \, \gamma_1\, i + \gamma_2,\quad i=1,2,\ldots\,. 
\en}

For the case $I=\{1,\dots,n\}$, the main result of \cite{BP} implies that the system \eqref{sdeI} has a weak solution (``weak existence"), which is unique in the sense of the probability distribution (``weak uniqueness"); the strict positivity of the diffusion coefficients is crucial here. In the case   $\,I=\nn\,$, a slight modification of the proof of Proposition 3.1 in \cite{Sh} shows that weak existence and weak  uniqueness hold  for the system \eqref{sdeI} under the Assumption \ref{asmp}; see Proposition \ref{inf_syst} below. 

\smallskip
In this paper, we investigate the questions of  existence of a strong solution (``strong existence") and  of pathwise uniqueness (``strong uniqueness") in both cases. Due to the discontinuity of the diffusion functions 
$$
\sigmab_i (x) \,=\, \sum_{j\in I} \mathbf{1}_{\{x_i =x_{(j)} \}}\, \sigma_j\,, \qquad x \in \mathbb{R}^n\,, \quad i=1, \ldots, n 
$$
in (\ref{sdeI}), general results on strong existence and strong uniqueness, which rely on the regularity of the diffusion coefficients, do not apply even when $I$ is finite.

\medskip

In order to construct a strong solution to the system \eqref{sdeI} in the case  $\,I=\{1,\ldots,n\}\,$, we rely heavily on the results of the recent  article \cite{IFKP}; this paper deals with the case $\,n=2\,$ and establishes   strength and pathwise uniqueness for the solution of the resulting system (\ref{sdeI}) (actually, even   when one of the diffusion coefficients vanishes, but not both).  The idea, then, is to put together paths of the strong solutions found in \cite{IFKP} for two particles, to obtain the strong solution in the case $n>2$ of several particles. 

This is   possible and the approach is viable, however, only  when the particle system in \eqref{sdeI} does not  exhibit triple collisions, that is, when the event
\eq\label{triple_coll_event}
\{\exists~ t>0\,:\;X_i(t)=X_j(t)=X_k(t)\;\,\text{for\;some\;}i<j<k\,\}
\en  
has zero probability for the state process $\,X\,$ in the weak solution of the system \eqref{sdeI}. We provide new, necessary and sufficient conditions for the absence of triple collisions in the case $\,I=\{1,\ldots,n\}\,$;  and develop the first such conditions in the case $\,I=\nn\,$. 

To formulate our main results we shall need the following Definition \ref{conc_seq}, as well as Conditions \ref{cond1} and \ref{cond2} below.

\begin{defn}
  \label{conc_seq}
A finite or infinite sequence $(a_1,a_2,\ldots)$ is called \textit{concave}, if for every three consecutive elements $a_i$, $a_{i+1}$, $a_{i+2}$   we have $$a_{i+1}\geq\frac{1}{\,2\,}\, \big(a_i+a_{i+2}\big)\,.$$ 
\end{defn} 

\begin{cond}\label{cond1}
{\it The sequence $(\sigma_1^2,\sigma_2^2,\ldots)$ is concave.}
\end{cond}

\begin{cond}
 \label{cond2}
{\it Either 
$\,I=\{1,2,\ldots,n\}\,$ and the sequence $\,(0,\sigma_1^2,\sigma_2^2,\ldots,$ $\sigma_n^2,0)\,$ is concave; or $\,I=\nn\,$ and the sequence $(0,\sigma_1^2,\sigma_2^2,\ldots)$ is concave. }
\end{cond}

Our main results read as follows. 

{\thm\label{triple_coll_thm} Consider the particle system in \eqref{sdeI} and, if $I=\nn$, let Assumption \ref{asmp} be satisfied. 

If the diffusion coefficients satisfy Condition \ref{cond2}, then the unique weak solution of \eqref{sdeI} has no triple collisions; that is, the event \eqref{triple_coll_event} has zero probability. On the other hand, if Condition \ref{cond1} fails, then the event \eqref{triple_coll_event} has   positive probability.}

{\thm\label{strong_sol} Consider the particle system in \eqref{sdeI} and, if $I=\nn$, let Assumption \ref{asmp} be satisfied. Introduce the first time of a triple collision, namely 
\eq\label{whatistau}
\tau \, := \, \inf\{ \, t\geq0 \,|\;\exists ~ i<j<k:\;X_i(t)=X_j(t)=X_k(t)\}\,.
\en
Then the system \eqref{sdeI} has a unique strong solution up to time $\tau$. 

In particular, if Condition \ref{cond2} is also satisfied, then there is a unique strong solution of the system \eqref{sdeI} defined for all $\,t\geq0$.}

\bigskip

\begin{rmk}
Theorem \ref{triple_coll_thm} can be recast as saying that Condition \ref{cond1} is necessary, and Condition \ref{cond2}  sufficient,  for the absence of triple collisions. A condition that is both necessary and sufficient  for the absence of triple collisions, has yet to be determined. So far, this question is completely resolved only in the case $n=3$, in which the results of \textsc{Varadhan \& Williams} \cite{VW} imply that Condition \ref{cond1} is both necessary and sufficient  for the absence of triple collisions; see the proofs of Lemma \ref{coll_all} and Theorem \ref{triple_coll_thm} below for more details. 

\medskip

The remaining gap between Condition \ref{cond1} and Condition \ref{cond2} is due to the following reason. By an inductive argument, we reduce the statement of Theorem \ref{triple_coll_thm} to the problem studied in \textsc{De Blassie} \cite{deB}. However, the (sharp) criterion given there involves the invariant distribution of the projection of a certain diffusion process in a Euclidean space on the unit sphere. Due to the lack of rotational symmetries in our situation, it is however not clear how to analyze this invariant distribution. For this reason, we simplify the condition in \cite{deB} to a checkable sufficient condition in Proposition \ref{deB_thm} below, sacrificing its sharpness at this point. 

\medskip

Theorem \ref{strong_sol} leaves open the questions  of whether a strong solution continues to exist, and of whether pathwise  uniqueness continues to hold, after a triple collision (we know from the work of \textsc{Bass \& Pardoux} \cite{BP} that a weak solution exists after such triple collisions, and is unique in distribution). At the moment, we conjecture that strong solutions fail to exist beyond the time of the first triple collision, but this problem remains open and will have to be settled in future work. $~~~~~~\Box$ 
\end{rmk}

\medskip
Questions regarding the presence or absence of triple or higher-order collisions in multidimensional diffusions have been addressed by several authors; in addition to the example of section 3 in \textsc{Bass \& Pardoux} \cite{BP}, one should consult the works by \textsc{Friedman} \cite{Fr} (see also \cite{Fr2}, chapter 11), \textsc{Ramasubramanian} \cite{Ram}, \cite{Ram2}, \textsc{De Blassie} \cite{deB}, \cite{deB2} and \textsc{Ichiba \& Karatzas} \cite{IK}. In the papers \cite{CL1} and \cite{CL2}, \textsc{C\'epa \& L\'epingle} consider systems of Brownian particles with repulsive forces of electrostatic type, and show absence of triple collisions under conditions of sufficiently strong repulsion. 

\bigskip
\noindent
{\it Preview:} The present paper is organized as follows. Section \ref{sec2} collects a number of preliminaries, most notably results from \cite{BFK}, \cite{IPBFK}, \cite{deB} and \cite{Sh} that are crucial in our context. Section \ref{sec3} deals with the absence of triple collisions under Condition \ref{cond2}, and with the proof of Theorem \ref{triple_coll_thm}, which represents a significant improvement over our earlier results in \cite{IK}. In particular, we provide here new  necessary conditions and new sufficient conditions for the absence of triple collisions in the case of  a finite number of paticles,  and develop the first such conditions for a countable infinity of particles. 

Section \ref{sec4} is devoted to the proof of Theorem \ref{strong_sol}. We start by setting up an inductive procedure, which ``bootstraps" the strength of the solution to the system of equations (\ref{sdeI}) that was  established recently by \textsc{Fernholz, Ichiba, Karatzas \& Prokaj}    \cite{IFKP} for the case $\, n=2\,$ of two particles --  first to the case $\, n=3\,$ of three particles;   then to the case of an arbitrary, finite number $\,n\,$ of particles;  and finally, building on results of \textsc{Shkolnikov} \cite{Sh}, to the case of a countable infinity of particles.   

\bigskip

\section{Preliminaries}
 \label{sec2}

We start with some preliminaries on the weak solution of the system \eqref{sdeI}, when $I=\{1,\dots,n\}$. First, we recall the dynamics of the ordered particles $X_{(1)},\dots,X_{(n)}$ in the system \eqref{sdeI} from section 3 in \cite{BFK} and section 4 of \cite{IPBFK}; once again, the strict positivity of the diffusion coefficients is crucial for these results. 

As in those papers, we shall denote by  $\,\Lambda^{j-1,j}(t)$, $j=2,\ldots,n\,$   the local times (normalized according to Tanaka's formula) accumulated at the origin by the nonnegative semimartingales 
\eq
 \label{spacey}
Y_{j-1} (\cdot)\,:=\, X_{(j)}(\cdot)-X_{(j-1)}(\cdot)\,, \quad j=2,\dots,n
\en
 over the inteval $\, [0,t]\,$, and   set $\,\Lambda^{0,1}(\cdot)\equiv \Lambda^{n,n+1}(\cdot)\equiv 0\,$.

{\prop\label{ordered_part} Set $I=\{1,\dots,n\}$ and let $(X,W)$ be a weak solution of the system \eqref{sdeI}. Then there exist independent standard Brownian motions $\,\beta_1,\ldots,\beta_n\,$ such that 
\eq
 \label{BG}
\mathrm{d}X_{(j)}(t)=\delta_j \;\mathrm{d}t + \sigma_j\;\mathrm{d}\beta_j(t) 
+ \frac{1}{\,2\,}\Big(\mathrm{d}\Lambda^{j-1,j}(t) - \mathrm{d}\Lambda^{j,j+1}(t)\Big),\;~~t\geq0  
\en
holds for all $j\in I\,$. 
}

\medskip

It was observed in \cite{BFK} that Proposition \ref{ordered_part} permits the identification of the process of ordered particles $(X_{(1)},\dots,X_{(n)})$ as a Reflected Brownian Motion (RBM for short) in the wedge
\eq
\mathfrak{W}:=\{(x_1,\dots,x_n)\in\rr^n:\;x_1\leq\ldots\leq x_n\}
\en
with reflection matrix
\begin{eqnarray*}
\mathfrak{R}\,:=\,\left(\begin{array}{ccc}
-\frac{1}{2} & 0 & 0 \\
\frac{1}{2} & -\frac{1}{2} & 0 \\
0 & \frac{1}{2} & \ddots \\
0 & 0 & \ddots
\end{array}\right).
\end{eqnarray*}
That is, the process $(X_{(1)},\dots,X_{(n)})$ behaves as an $n$-dimensional standard Brownian motion in the interior of the wedge $\,\mathfrak{W}\,$, and is obliquely reflected on the faces $\{x_i=x_{i+1}\}$, $i=1,\ldots,n-1$ of $\,\mathfrak{W}\,$. The directions of reflection are specified by the columns $\mathfrak{r}_i$, $i=1,\ldots,n-1$ of the matrix $\mathfrak{R}\,$.   

\medskip

Occasionally, it will be more convenient to consider instead of the process of the ordered particles $(X_{(1)},\ldots,X_{(n)})$ the {\it process of spacings} (gaps) 
\eq
Y\,:=\, \big(X_{(2)}-X_{(1)},\ldots,X_{(n)}-X_{(n-1)}\big) 
\en
as in (\ref{spacey}). From Proposition \ref{ordered_part}, we have the dynamics 
\eq
 \label{spacey_dyn}
\mathrm{d}Y_{j-1}(t)\,=\,\big( \delta_j - \delta_{j-1} \big) \;\mathrm{d}t + \sigma_j\;\mathrm{d}\beta_j(t) - \sigma_{j-1}\;\mathrm{d}\beta_{j-1}(t) \qquad  \qquad 
\en
$$
\qquad  \qquad \qquad - \;\frac{1}{\,2\,}\Big(\mathrm{d}\Lambda^{j ,j+1}(t) + \mathrm{d}\Lambda^{j-2,j-1}(t)\Big)+  \mathrm{d}\Lambda^{j-1,j}(t)\,,\;\quad t \ge 0
$$
for the spacings of (\ref{spacey}) with $\,j=2, \ldots, n\,$. Thus, 
the process $Y$ is an RBM   in the $(n-1)$-dimensional orthant $(\rr_+)^{n-1}$ with reflection matrix
\begin{eqnarray*}
\mathcal{R}\,:=\,\left(\begin{array}{cccc}
1 & -\frac{1}{2} & 0 & 0 \\
-\frac{1}{2} & 1 & -\frac{1}{2} & 0 \\
0 & -\frac{1}{2} & 1 & \ddots \\
0 & 0 & \ddots & \ddots
\end{array}\right).
\end{eqnarray*}
For a detailed summary of many results on Brownian motions with oblique reflection in the orthant, we refer to the excellent survey article \cite{W1}.

\medskip

For further reference we make the following simple observation. The event in \eqref{triple_coll_event} can be reformulated as
\eq
\{\,\exists~ t\geq0:\;Y_i(t)=Y_{i+1}(t)=0\;\text{for\;some}\;1\leq i\leq n-1\}\,,
\en
so that the presence or absence of triple collisions is an intrinsic property of the spacings process $Y$.

\bigskip

Next, we let $\,I=\nn\,$ and construct the weak solution to \eqref{sdeI} along the lines of the proof of Proposition 3.1 in \cite{Sh}, under the Assumption \ref{asmp}. 

{\prop\label{inf_syst} Let $I=\nn$ and let Assumption \ref{asmp} be satisfied. 

There exists then a  weak solution $(X,W)$ of the system \eqref{sdeI}, which is unique in distribution. 

Moreover, after enlarging the probability space if necessary, we can find stopping times $0=\kappa_0\leq\kappa_1\leq\ldots$, integers $M=n(0)<n(1)<\ldots$ and weak solutions $X^{(k)}$, $k\geq0$ of the system \eqref{sdeI} with $I=\{1,\dots,n(k)\}$ such that
\begin{eqnarray}
&&X_i(t)=X^{(k)}_i(t),\quad t\in[\kappa_k,\kappa_{k+1}],\;1\leq i\leq n(k),\\
&&X_i(t)=X_i(0)+\delta_i t + \sigma_i W_i(t),\quad t\in[\kappa_k,\kappa_{k+1}],\;i\geq n(k)+1
\end{eqnarray}
and, for each $\,k\geq0\,$, the processes $X^{(k)}$ and $(W_{n(k)+1},W_{n(k)+2},\dots)$ are independent.}

\bigskip

\noindent{\bf Proof.} The proof of Proposition 3.1 in \cite{Sh} carries over {\it mutatis mutandis} to the situation here. We only need to replace the a priori estimate on the expected number of particles in an interval of the form $(-\infty,x]$ at a time $\,t\geq0\,$ by  
\[
\sum_{i \in \nn} \,\sup_{\varsigma(\cdot)} \,
\mathbb{P}\left(X_i(0)-\max_i |\delta_i|\cdot t - \sup_{0\leq s\leq t} \int_0^s \varsigma(u)\;\mathrm{d}W_i(u)<x\right)<\infty\,,
\]
where the supremum is taken over all progressively measurable processes $\,\varsigma(\cdot)\,$ adapted to the filtration on the underlying probability space, which take values in the interval $\,[\,\min_i \sigma_i\,,\, \max_i \sigma_i\,]\,$. The finiteness of the series follows from the $\,\mathbb{L}^2$-version of Doob's maximal inequality for continuous martingales and the condition \eqref{lin_growth}. \ep 

\bigskip

Finally, as a preparation for the proof of Theorem \ref{triple_coll_thm}, we state a special case of the main result of \textsc{De Blassie} in \cite{deB}. 

{\prop\label{deB_thm} Let $\,\mathbb{H}$ be a finite-dimensional Euclidean space and $\,\sigmab: \mathbb{H}\rightarrow \mathbb{H}^2$ be a mapping satisfying $$\sigmab(x)\,=\,\sigmab \left(\frac{x}{\,\,\|x\|_2\,}\right) \qquad \hbox{for all} \quad x\in \mathbb{H} \setminus \{0\}\, $$ where $\, \|\cdot\|_2$ denotes the Euclidean norm on $\,\mathbb{H}\,$, and suppose that the set of discontinuity points of $\,\sigmab (\cdot)\,$ on the unit sphere $\{x\in H:\,\|x\|_2=1\}$ has surface measure zero. In addition, let $\,W$ be an $\,\mathbb{H}-$valued standard Brownian motion, and suppose that the martingale problem corresponding to the $\,\mathbb{H}-$valued stochastic differential equation 
\eq
\mathrm{d}V(t)\,=\,\sigmab \big(V(t)\big)\;\mathrm{d}W(t)\,, \qquad V(0) \in \mathbb{H} \setminus \{ 0\}
\en
is well-posed. If the condition
\eq\label{deB_cond}
\inf_{x\in \mathbb{H} \atop \|x\|_2=1} \left( \frac{\mathrm{tr}\;\alphab (x)}{\langle \alphab (x)x,x\rangle} \right)\, >\, 2
\en
is satisfied, then we have 
$$
\mathbb{P} \, \big( \,V(t) \neq 0\,, \quad \forall \; t \in [0, \infty)\, \big)\,=\,1\,.
$$
Here  $\,\alphab (\cdot)=\sigmab(\cdot)'\sigmab(\cdot)$ is the diffusion matrix of $\,V$, $\,\mathrm{tr}$ denotes the trace operator, and $\langle \cdot \,,\cdot \rangle$ is the Euclidean scalar product on $\,\mathbb{H}\,$.} 

\bigskip

\noindent{\bf Proof.} It suffices to note that (\ref{deB_cond}) implies the condition $$\inf_{x\in \mathbb{H} \atop \|x\|_2=1} B(x)>1$$ in the notation of equation (1.9) in \cite{deB}. Thus, the result is a special case of Theorem 1.1 (i) in \cite{deB}. \ep

\medskip

\section{Triple collisions}
 \label{sec3}

The two main steps in the proof of Theorem \ref{triple_coll_thm} are provided by the following two lemmas. 

{\lemma\label{coll_all} Let $\,I=\{1,\dots,n\}\,$ with an integer $\,n\geq3\,$, and suppose that Condition \ref{cond2} holds. Then the first time of a triple collision $\tau$, defined in \eqref{whatistau}, must satisfy
\eq
\tau\,=\,\eta
\en
with probability one, where 
\eq
 \label{eta}
\eta\,:=\,\inf\{t\geq0:\;X_1(t)=X_2(t)=\ldots=X_n(t)\}\,.
\en
}

{\lemma\label{no_triple_coll} Let $I=\{1,\dots,n\}$ with an $n\geq3$ and suppose that Condition \ref{cond2} holds. Then the first time of a triple collision $\tau$, defined in \eqref{whatistau}, must satisfy $\,\tau=\infty\,$ with probability one.}

\bigskip

We shall prove both lemmas simultaneously, by induction over $n\,$.

\bigskip

\noindent{\bf Proof of Lemmas \ref{coll_all} and \ref{no_triple_coll}:} {\bf Step A.} First, we note that for every $\,T\in (0, \infty)\,$ we can make a Girsanov change of measure such that the processes ($\delta_i \, t + \sigma_i \, \beta_i(t)$, $t\in[0,T]$), $i=1,\ldots,n$ become independent standard Brownian motions under the new measure. Here, the Brownian motions $\beta_1,\ldots,\beta_n$ are defined as in Proposition \ref{ordered_part}. 

Since almost-sure  events   retain this property under an   absolutely continuous   change of  measure, and since $\,T\in (0, \infty)\,$ can be chosen arbitrarily, it suffices to prove the almost sure absence of triple collisions under the new measure. For notational simplicity, we prefer to assume $\delta_1=\ldots=\delta_n=0$ from the start. Section 2.2 in \cite{IK} can be consulted for a more detailed exposition of the same argument.   
 
\bigskip
\noindent{\bf Step B.} We proceed with the inductive argument. For $n=3$, we deduce from  Proposition 3 of \cite{IK} that
\eq
 \tau\,=\,\eta \,=\, \infty
\en
holds with probability one, in the notation of \eqref{whatistau} and (\ref{eta}). In fact, this is a consequence of Theorem 2.2 in \cite{VW} for the reflected Brownian motion $Y$ (see the proof of Proposition 3 in \cite{IK} for more details).  

\bigskip
\noindent{\bf Step C.} Now, fix an $m\geq4$ and suppose that Lemmas \ref{coll_all} and \ref{no_triple_coll} hold for all $3\leq n<m$. We will first show that Lemma \ref{coll_all} must hold for $n=m$ as well. To this end, we define for each $0<\varepsilon<1$ the stopping time 
\eq
 \label{eta_epsilon}
\eta_\varepsilon\,:=\,\inf\{t\geq0:\;\|Y(t)\|_2\leq\varepsilon\;\;\text{or}\;\;\|Y(t)\|_2\geq\varepsilon^{-1}\},
\en
where we have written $\|\cdot\|_2$ for the usual Euclidean norm. 

We claim that, for all $\, 0<\varepsilon<1\,$,   the comparison  
$$
\tau\geq\eta_\varepsilon~~\; \hbox{holds with probability one}  \,.
$$ 
If we can prove this claim, then we will be able to conclude that either $\,\tau=\eta=\infty\,$ or $\,\tau=\eta=\lim_{\varepsilon\downarrow0}\eta_\varepsilon<\infty\,$ must hold, and this  will yield Lemma \ref{coll_all}.

\medskip

To prove the claim, we deploy an argument similar to the one on pages 471-475 in \cite{W2}. As there, we consider the local behavior of the RBM of interest (in our case $Y$) on the compact sets
\eq
 \label{K}
K_\varepsilon\,:=\, \big\{ \, y\in(\rr_+)^{n-1}:\;\varepsilon\leq\|y\|_2\leq\varepsilon^{-1} \, \big\}\,,\quad 0<\varepsilon<1\, .
\en
For each $y\in(\rr_+)^{n-1} \setminus \{0\}\,$, we define an open set $U(y)$ of $(\rr_+)^{n-1}$ such that 
\eq\label{nbhd}
\exists ~1\leq j\leq n-1,\;\delta>0:\;z_j\geq\delta\;\text{for\;all}\;\;z\in U(y)\,,
\en
and let $\,\mathrm{j}(y)\,$ be a number in $\{1,\dots,n-1\}$ as in \eqref{nbhd}. From the semimartingale decomposition of $\,Y$ in (\ref{spacey_dyn}), we see that if we start $Y$ in the set $U(y)$ for some $\, y\in(\rr_+)^{n-1} \setminus \{0\}\,$, then we can write
\begin{eqnarray*} 
&&\big(Y_j(t\wedge\zeta_{U(y)})\,:\,1\leq j\leq \mathrm{a}(y)\big)=\big(Y'_j(t\wedge\zeta_{U(y)}):\,1\leq j\leq \mathrm{a}(y)\big),\\
&&\big(Y_j(t\wedge\zeta_{U(y)}):\,\mathrm{b}(y)\leq j \leq n-1\big)=\big(Y''_j (t\wedge\zeta_{U(y)}):\,1\leq j \leq \mathrm{c}(y)\big)
\end{eqnarray*}
for all $\,t\geq0\,$. Here we have set $$\mathrm{a}(y):=\mathrm{j}(y)-1\,, \quad \mathrm{b}(y):=\mathrm{j}(y)+1\,, \quad \mathrm{c}(y):=n-1-\mathrm{j}(y)\,;$$ the process $\,Y'\,$ is an RBM in the $\,\mathrm{a}(y)-$dimensional orthant; the process $\,Y''\,$ is an RBM in the $\,\mathrm{c}(y)-$dimensional orthant; and $\zeta_{U(y)}$ is the time that $Y$ hits the boundary of $U(y)$. In particular, the induction hypothesis implies 
\eq
\tau>\zeta_{U(y)}\,,
\en
where $y$ is such that $Y(0)\in U(y)$.

\bigskip

Next, we fix an $\varepsilon\in(0,1)$ and cover the compact set $K_\varepsilon$ of (\ref{K}) by a finite number of open sets from the collection $\,U(y)\,$, $\,y\in(\rr_+)^{n-1} \setminus \{0\}\,$, say
\eq
K_\varepsilon \, \subset \, \bigcup_{\ell=1}^L \,U(y_\ell)\,. 
\en
Then we can find stopping times $\, 0= \zeta_0<\zeta_1<\zeta_2<\ldots\,$ of the form $\zeta_{U(y)}$, such that the path of the process $Y(t\wedge\eta_\epsilon)$, $t\geq0$ can be decomposed into 
\eq
Y(t\wedge\eta_\epsilon)\,,\;\;\zeta_k\leq t\leq\zeta_{k+1}\,,
\en
$k\geq0$ with the notation of (\ref{eta_epsilon}) and with $\,Y(t\wedge\eta_\varepsilon)\in\overline{U(y_\ell)}$ for some $1\leq \ell\leq L$ and all $\zeta_k\leq t\leq\zeta_{k+1}$. 

\smallskip
Using the strong Markov property of $Y$ and the previous observation, one shows by induction over $k$ that $\tau>\zeta_k\wedge\eta_\varepsilon$ must hold with probability one, for all $k\geq0$. By taking the limit $k\rightarrow\infty$, we conclude that $\,\tau\geq\eta_\varepsilon\,$ holds  with probability one. Thus, outside of a set of probability zero, we must have 
\eq
\forall \; \varepsilon\in(0,1):\;\tau\geq\eta_\varepsilon.
\en
This proves the claim and, thus, Lemma \ref{coll_all} for $n=m$. 

\bigskip 

\noindent{\bf Step D.} It remains to show that Lemma \ref{no_triple_coll} holds for $n=m$. To this end, consider the centered process
\eq
V(t)\,:=\, \bigg(X_1(t)-n^{-1}\sum_{i=1}^n X_i(t),\cdots,X_n(t)-n^{-1}\sum_{i=1}^n X_i(t)\bigg)
\en
for $\, 0 \le t < \infty\,$. It is obvious that
\eq
\eta\,=\,\inf\{t\geq0:\;V(t)=0\}.
\en
In addition, recalling that without loss of generality we have assumed $\delta_1=\ldots=\delta_n=0$ in Step A, we see from \eqref{sdeI} that $V$ is a diffusion process in the hyperplane
\eq
\mathbb{H}\,=\,\{\,x\in\rr^n:\;x_1+\ldots+x_n=0\,\}
\en
with zero drift, and with diffusion function satisfying the conditions of Proposition \ref{deB_thm}. Moreover, the martingale problem corresponding to the stochastic differential equation for the process $V$ is well-posed, thanks to the main result of \cite{BP} and the non-degeneracy of the diffusion matrix of $V$ on $\mathbb{H}$ (which, in turn, follows from the strict positivity of the diffusion coefficients in \eqref{sdeI}). Thus, we  only need to check that condition \eqref{deB_cond}  is a consequence of  Condition \ref{cond2}, since then Lemma \ref{no_triple_coll} will follow from Proposition \ref{deB_thm}.   

\bigskip

To check \eqref{deB_cond}, we first compute the diagonal entries of the diffusion matrix $\, \alphab (\cdot)\,$ of $V$ in the coordinates of $\,\rr^n\,$ to
\eq
\sigma_i^2 (1-2n^{-1}) + n^{-2}\sum_{j=1}^n \sigma_j^2\,\,,\quad i=1,\ldots,n\,,
\en
where  the order depends on the ranking of the coordinates of $\,V$. Next, we note that the normal vector to $\mathbb{H}$ is an eigenvector of $\, \alphab(\cdot)\,$ with eigenvalue $0\,$. Thus, $\, \alphab(\cdot)\,$ has an orthonormal eigenbasis over $\rr^n$, which includes the unit normal vector to $\mathbb{H}\,$. It follows that the trace of the diffusion matrix $\, \alphab(\cdot)\,$ of $V$ on $\,\mathbb{H}\,$ coincides with the trace of the diffusion matrix of $V$ in the coordinates of $\,\rr^n\,$, and is given by      
\[
M_0\,:=\,\sum_{i=1}^n \Big(\sigma_i^2 (1-2n^{-1}) + n^{-2}\sum_{j=1}^n \sigma_j^2\Big)\,=\,\frac{n-1}{n}\sum_{i=1}^n \sigma_i^2\,.
\]

Now, we estimate the denominator on the left-hand side of \eqref{deB_cond} from above by the maximal eigenvalue of the diffusion matrix $a(\cdot)$ of $\,V$ on $\,\mathbb{H}\,$. The latter is given by the maximal eigenvalue of the matrix 
\eq
 \label{PP}
P^{\,\prime} \, \mathrm{diag}(\sigma_1^2,\ldots,\sigma_n^2)\,P\,,
\en
 where $P$ is the matrix representing the normal projection of vectors in $\rr^n$ onto $\,\mathbb{H}\,$, and $\mathrm{diag}(\sigma_1^2,\ldots,\sigma_n^2)$ is the $n\times n$ diagonal matrix with diagonal entries $\sigma_1^2,\ldots,\sigma_n^2$. 
 
 \smallskip
 Next, we observe that the spectral radius (and, thus, the maximal eigenvalue) of the matrix (\ref{PP})  is given by  
\eq
 \label{maximum1}
M_1:=\max_{x\in\rr^n \setminus \{0\} \atop x_1+\ldots+x_n=0} \left( \frac{\sigma_1^2 x_1^2+\ldots+\sigma_n^2 x_n^2}{x_1^2+\ldots+x_n^2} \right)\,. 
\en
This quantity can be further estimated from above as  
\eq
 \label{maximum2}
M_2\,:=\, \max_{x\in\rr^n \setminus \{0\} \atop x_1+\ldots+x_n=0} \left( \frac{C x_1^2+ c (x_2^2 + \ldots + x_n^2)}{x_1^2+\ldots+x_n^2} \right)\, \ge \, M_1\,,
\en
where $C$ is the maximal element of the set $\{\sigma_1^2,\ldots,\sigma_n^2\}$ and $c$ is the second largest element of the same set. A careful optimization using Lagrange multipliers gives $$M_2\,=\,\frac{\,n-1\,}{n}\,C+\frac{1}{\,n\,}\, c\,,$$ and shows that $M_2=M_1$ holds if and only if all elements in the set $\{\sigma_1^2,\ldots,\sigma_n^2\}$ are greater than or equal to $c$. 

\smallskip
All in all, we conclude that the condition \eqref{deB_cond}, which proscribes triple collisions,  amounts to 
\eq
 \label{deB_cond3}
M_0 > 2\,M_1
\en
and is satisfied, in particular, if the stronger inequality
\eq
 \label{deB_cond2}
M_0=\frac{n-1}{n}\sum_{i=1}^n \sigma_i^2 \,>\, 2\, \Big(\,\frac{n-1}{n}\,C+\frac{1}{n}\,c\Big)= 2\,M_2
\en
holds. We shall show \eqref{deB_cond3}, or its stronger version \eqref{deB_cond2}, by distinguishing the cases $\sigma_1^2\neq C\neq \sigma_n^2$ and $C=\sigma_1^2$ (the case $C=\sigma_n^2$ being completely analogous to the latter). 

\smallskip
\noindent
$\bullet~$ 
In the first case ($\sigma_1^2\neq C\neq \sigma_n^2$), the concavity of the sequence $(0,\sigma_1^2,\ldots,\sigma_n^2,0)$ shows that  
its third-largest element is greater than or equal to $\frac{1}{\,2\,}\,C$, whereas its fourth-largest element is greater than or equal to $\frac{1}{\,3\,}\,C$. Hence, the left-hand side of \eqref{deB_cond2} is greater than or equal to 
\[
\frac{n-1}{n}\Big(C+c+\frac{1}{\,2\,}\,C+\frac{1}{\,3\,}\,C\Big)\,.
\]
Plugging this into \eqref{deB_cond2} and simplifying, we see that it suffices to check $\,(n-3)\,c\,>\,\frac{\,n-1\,}{6}\,C\,$. The latter inequality, thus also (\ref{deB_cond2}) and (\ref{deB_cond}), holds for all $\,n\geq4\,$ due to $\,c\geq\frac{2}{3}C\,$, a consequence of the concavity of the sequence $(0,\sigma_1^2,\ldots,\sigma_n^2,0)$. 

\smallskip
\noindent
$\bullet~$ 
In the second case ($C=\sigma_1^2$), we use the concavity of the sequence $\,(0,\sigma_1^2,\ldots,\sigma_n^2,0)$ to deduce that  
its third-largest element is greater than or equal to $\frac{1}{\,2\,}\,C$, whereas its fourth-largest element is greater than or equal to $\frac{1}{\,4\,}\,C$. Hence, the left-hand side of \eqref{deB_cond2} is greater than or equal to 
\[
\frac{n-1}{n}\,\Big(C+c+\frac{1}{\,2\,}\,C+\frac{1}{\,4\,}\,C\Big)\,.
\]   
Plugging this into \eqref{deB_cond2}, we conclude that it suffices to show  
\eq
 \label{strict}
(n-3)\,c\,>\,\frac{\,n-1\,}{4}\,C \,.
\en 
Using $\,c\geq\frac{3}{\,4\,}\,C\,$ (again, a consequence of the concavity of $(0,\sigma_1^2,\ldots,\sigma_n^2,0)\,$), we observe:
\eq
 \label{nonstrict}
(n-3)\,c\,\ge\,\frac{\,n-1\,}{4}\,C\,,\qquad \forall \,\; n \ge 4\,,
\en 
with equality if and only if both $\,c=\frac{3}{\,4\,}\,C\,$ and $\,n=4\,$ hold. Moreover, the derivation of \eqref{nonstrict} shows that: either \eqref{deB_cond2} holds (therefore, also \eqref{deB_cond3}); or we have $\,n=4\,$, $\,(\sigma_1^2,\sigma_2^2,\sigma_3^2,\sigma_4^2)=(C,\frac{3}{\,4\,}\,C,\frac{1}{\,2\,}\,C,\frac{1}{\,4\,}\,C)\,$ and equality in \eqref{deB_cond2}. In this latter case, however, the inequality $M_2 \geq M_1$ in (\ref{maximum2}) is strict, so that \eqref{deB_cond3} must hold. 
\ep 

\bigskip

We can now combine our results, to prove Theorem \ref{triple_coll_thm}. 

\bigskip

\noindent{\bf Proof of Theorem \ref{triple_coll_thm}:} {\bf Step 1.} First, let $I=\{1,\dots,n\}$. Then under Condition \ref{cond2} there are no triple collisions,  by virtue of Lemma \ref{no_triple_coll}. 

Now, suppose that Condition \ref{cond1} fails; that is, for some integer $\, i= 2,\ldots,n-1\,$  the comparison 
\eq
\label{fail}
 \sigma_i^2-\sigma_{i-1}^2 < \sigma_{i+1}^2-\sigma_i^2 
 \en
  holds. Consider the weak solution $X'=(X'_1,X'_2,X'_3)$ to the system \eqref{sdeI} with $I=\{1,2,3\}$ and the parameters $\,\delta_{i-1},\delta_i,\delta_{i+1},\sigma_{i-1},\sigma_i,\sigma_{i+1}\,$. Then, applying Theorem 2.2 of \cite{VW} as in the proof of Proposition 3 in \cite{IK}, we conclude that a triple collision of the particles $\,X'_1,\,X'_2,\,X'_3\, $ occurs with positive probability. It follows that there is a $T\in (0, \infty)\,$ and a bounded, open subset $U$ of the wedge $\,\{x\in\rr^3:\;x_1\leq x_2\leq x_3\}\,$, such that the event
\begin{eqnarray*}
&&\big\{ \, (X'_{(1)}(t),X'_{(2)}(t),X'_{(3)}(t))\in U,~t\in[0,T]\, \big\}\\
&\cap&\big\{\, \exists ~ t\in[0,T]:\,X'_1(t)=X'_2(t)=X'_3(t)\, \big\} 
\end{eqnarray*} 
has   positive probability for every initial condition in $U$. Along with the semimartingale decomposition of (\ref{BG}) for the components of the process $\,(X_{(1)},\dots,X_{(n)})\,$, this  implies that the event 
\begin{eqnarray*}
&&\big\{ \,X_{(j-2)}(t)\neq X_{(j-1)}(t),X_{(j+1)}(t)\neq X_{(i+2)}(t),~t\in[0,2T]\,\big\}\\
&\cap&\big\{\, (X_{(j-1)}(t),X_{(j)}(t),X_{(j+1)}(t))\in U,~t\in[T,2T]\, \big\}\\
&\cap&\big\{\,\exists ~ t\in[T,2T]:\,X_{(j-1)}(t)=X_{(j)}(t)=X_{(j+1)}(t)\, \big\}
\end{eqnarray*}
has  positive probability. This completes the proof of Theorem \ref{triple_coll_thm} for $I=\{1,\dots,n\}$.

\bigskip

\noindent{\bf Step 2.} Now, we turn to the case $\,I=\nn\,$ and assume first that Condition \ref{cond2} holds. We recall the notation  in Proposition \ref{inf_syst} and observe that the event \eqref{triple_coll_event} is contained in the event
\eq
\bigcup_{i_1<i_2<i_3} \, \bigcup_{k \in \mathbb{N}_0}\, \big\{ \,X_{i_1}(t)=X_{i_2}(t)=X_{i_3}(t)\;\text{for\;some\;}t\in[\kappa_k,\kappa_{k+1}] \, \big\}\,.
\en
Moreover, for every fixed $\, i_1<i_2<i_3\,$ and $\,k\geq0\,$, Proposition \ref{inf_syst} shows that there is a choice of $n\geq M$ and a weak solution $\big(X^{(k)}_1,\cdots,X^{(k)}_n \big)$ of the system \eqref{sdeI} with $I=\{1,\dots,n\}$ and parameters $\delta_1,\dots,\delta_n,\sigma_1,\dots,\sigma_n\,$, such that 
\eq
X_i(t)=X^{(k)}_i(t),\quad t\in[0,\kappa_{k+1}],\;1\leq i\leq n.
\en  
Therefore, Step 1 of the present proof implies 
\eq
\mathbb{P}\Big(\big\{ \, X_{i_1}(t)=X_{i_2}(t)=X_{i_3}(t)\;\;\text{for\;some\;\;}t\in[\kappa_k,\kappa_{k+1}] \, \big\}\Big)\,=\,0
\en
for any fixed $i_1<i_2<i_3$ and $k\geq0$. Thus, the event \eqref{triple_coll_event} has zero probability, as claimed.

\medskip
Next, suppose that Condition \ref{cond1} fails, and recall the definition of the constant $M$ in Assumption \ref{asmp}. In this case, there is an integer $\,i =2,\ldots,M\,$ such that \eqref{fail} holds. In particular, it follows from Step 1 of this proof that the weak solution $X'$ of \eqref{sdeI} with $I=\{1,2,\ldots,M+1\}$ and parameters $\delta_1,\delta_2,\ldots,\delta_{M+1}$, $\sigma_1,\sigma_2,\ldots,\sigma_{M+1}$ exhibits a   triple collision with positive probability. In conjunction with the semimartingale decomposition of the process $\,X'$, this shows that for every $L>0$ there is a $T>0$ such that the event
\begin{eqnarray*}
&&\{X'_{(M+1)}(t)\leq X'_{(M+1)}(0)-2L,\;T\leq t\leq 2T\}\\
&\cap&\{\exists ~i_1<i_2<i_3,\,t\in[T,2T]:\;X'_{i_1}(t)=X'_{i_2}(t)=X'_{i_3}(t)\}
\end{eqnarray*}
has a positive probability. In addition, Proposition \ref{inf_syst} shows that, with probability one, the paths of the process $X_{(M+2)}$ are continuous and do not reach negative infinity in finite time. 

Putting these last two observations together, we see that there are exist real constants $L>0$ and   $T>0\,$, such that the event
\begin{eqnarray*}
&& \big\{ \, \max(X_{(M+1)}(t),X_{(M+1)}(0)-L)<X_{(M+2)}(t),\;\;t\in[0,2T]\, \big\}\\
&\cap& \big\{ \,X_{(M+1)}(t)\leq X_{(M+1)}(0)-2L,\;T\leq t\leq 2T \, \big\}\\
&\cap&\big\{\,\exists ~ i_1<i_2<i_3\leq M+1,\,t\in[T,2T]:\;X_{i_1}(t)=X_{i_2}(t)=X_{i_3}(t)\,\big\}
\end{eqnarray*}
has a positive probability. In particular, the event \eqref{triple_coll_event} has a positive probability. \ep  

\section{Construction of strong solutions}
 \label{sec4}

This section is devoted to the proof of Theorem \ref{strong_sol}. In the first subsection we explain our methodology in the special case   $I=\{1,2,3\}$. The following subsection extends the construction of strong solutions to systems with any finite number of particles. Finally, in the last subsection we use the strong solutions in systems with finitely many particles to obtain the strong solution in the system with infinitely many particles.

\subsection{Systems with three particles}

In this subsection we explain the construction of strong solutions in the case that there are only three particles. That is, we consider the system of SDEs
\eq\label{sde3}
\mathrm{d}X_i(t)\,=\,\sum_{j=1}^3 \mathbf{1}_{\{X_i(t)=X_{(j)}(t)\}}\;\delta_j\;\mathrm{d}t 
\,+\, \sum_{j=1}^3 \mathbf{1}_{\{X_i(t)=X_{(j)}(t)\}}\;\sigma_j\;\mathrm{d} W_i(t)
\en
with initial conditions satisfying $X_1(0)<X_2(0)<X_3(0)$. As before, we define $\tau$ to be the first time of a triple collision for a weak solution of the system \eqref{sde3}. To wit,
\eq
\tau=\inf\{t\geq0:\;X_1(t)=X_2(t)=X_3(t)\}.
\en
In this setting our main result reads as follows. 

{\prop\label{strong3}
Suppose $$\sigma_2^2-\sigma_1^2\,\geq\,\sigma_3^2-\sigma_2^2\,.$$ Then the system \eqref{sde3} admits a strong solution, which is pathwise unique. Moreover, if $$\sigma_2^2-\sigma_1^2\,<\,\sigma_3^2-\sigma_2^2$$ holds, then a strong solution exists and is pathwise unique up to the triple collision time $\tau$.}

\bigskip

\noindent\textbf{Proof.} We will construct a strong solution to \eqref{sde3} by putting together paths of strong solutions to the system \eqref{sdeI} with $I=\{1,2\}$ and parameters $\,\delta_1,\delta_2,\sigma_1,\sigma_2\,$ and $\,\delta_2,\delta_3,\sigma_2,\sigma_3\,$, respectively. 

To this end, we introduce the following notation. For a closed time-interval $[a,b]$ we denote by $Z^{\,[a,b],B,W}(b_1,b_2,c_1,c_2)$ the $\rr^2$-valued strong solution of the system \eqref{sdeI} with drift parameters $b_1$, $b_2$ and diffusion  parameters $\,c_1$, $\,c_2\,$, driven by the independent standard Brownian motions $B$, $W$ on the time-interval $[a,b]$. The latter strong solution exists and is pathwise unique, on the strength of the results in section 5 of \cite{IFKP}. 

\bigskip

We can now define the stopping times $0=\tau_0\leq\rho_0\leq\tau_1\leq\rho_1\leq\ldots$ and the desired strong solution on the intervals $[\tau_k,\rho_k]$, $[\rho_k,\tau_{k+1}]$, $k\geq0$ inductively by 
\medskip
\begin{eqnarray*}
&&X_{\pi_k(1)}([\tau_k,\rho_k]):=\big(Z^{[\tau_k,\rho_{k }],W_{\pi_k(1)},W_{\pi_k(2)}}(\delta_1,\delta_2,\sigma_1,\sigma_2)\big)_1\,\,,\\
&&X_{\pi_k(2)}([\tau_k,\rho_k]):=\big(Z^{[\tau_k,\rho_{k }],W_{\pi_k(1)},W_{\pi_k(2)}}(\delta_1,\delta_2,\sigma_1,\sigma_2)\big)_2\,\,,\\
&&X_{\pi_k(3)}(t):=X_{\pi_k(3)}(\tau_k)+\delta_3 (t-\tau_k)+\sigma_3 W_{\pi_k(3)}(t) - \sigma_3 W_{\pi_k(3)}(\tau_k),\\
&&\quad\quad\quad\quad\quad\quad\quad\quad\quad\quad\quad\quad\quad\quad\quad\quad\quad\quad\quad\quad\quad\quad\;\; t\in[\tau_k,\rho_k]\\
&&\rho_k:=\inf\{t>\tau_k:\;X_{\pi_k(3)}(t)=X_{\pi_k(2)}(t)\;\text{or}\;X_{\pi_k(3)}(t)=X_{\pi_k(1)(t)}\}\,,
\end{eqnarray*}
\medskip
\begin{eqnarray*}
&&X_{\theta_k(1)}(t):=X_{\theta_k(1)}(\rho_k)+\delta_1 (t-\rho_k)+\sigma_1 W_{\theta_k(1)}(t) -\sigma_1 W_{\theta_k(1)}(\rho_k)\,,\\
&&\quad\quad\quad\quad\quad\quad\quad\quad\quad\quad\quad\quad\quad\quad\quad\quad\quad\quad\quad\quad\quad\quad t\in[\rho_k,\tau_{k+1}]\\
&&X_{\theta_k(2)}([\rho_k,\tau_{k+1}]):=\big(Z^{[\rho_k,\tau_{k+1}],W_{\theta_k(2)},W_{\theta_k(3)}}(\delta_2,\delta_3,\sigma_2,\sigma_3)\big)_1\,,\\
&&X_{\theta_k(3)}([\rho_k,\tau_{k+1}]):=\big(Z^{[\rho_k,\tau_{k+1}],W_{\theta_k(2)},W_{\theta_k(3)}}(\delta_2,\delta_3,\sigma_2,\sigma_3)\big)_2\,,\\
&&\tau_{k+1}:=\inf\{t>\rho_k:\;X_{\theta_k(2)}(t)=X_{\theta_k(1)}(t)\;\text{or}\;X_{\theta_k(3)}(t)=X_{\theta_k(1)}(t)\}.
\end{eqnarray*}

\smallskip
For each $\,k\geq0\,$, we have denoted here by $\,\pi_k\,$ a permutation of the set $\{1,2,3\}$ such that 
\[
X_{\pi_k(1)}(\tau_k)\leq X_{\pi_k(2)}(\tau_k)\leq X_{\pi_k(3)}(\tau_k)\,,
\]
and by $\, \theta_k\,$   a permutation of the set $\{1,2,3\}$ such that 
\[
X_{\theta_k(1)}(\rho_k)\leq X_{\theta_k(2)}(\rho_k)\leq X_{\theta_k(3)}(\rho_k)\,.
\] 

\medskip
It is straightforward to check that the just constructed processes $X_1,X_2,X_3$ are well-defined and form a strong solution to the system \eqref{sde3} up to the time 
\eq
\widetilde{\tau}:=\lim_{k\rightarrow\infty}\tau_k=\lim_{k\rightarrow\infty}\rho_k.
\en
On the other hand, since the paths of the ranked processes $X_{(1)}$, $X_{(2)}$, $X_{(3)}$ are continuous, we have
\[
X_{(1)}(\widetilde{\tau}\,)=\lim_{k\rightarrow\infty}X_{(1)}(\tau_k)=\lim_{k\rightarrow\infty} X_{(2)}(\tau_k)=X_{(2)}(\widetilde{\tau}\,). 
\]
Moreover, an analogous computation yields $\,X_{(2)}(\widetilde{\tau}\,)=X_{(3)}(\widetilde{\tau}\,)\,$. Thus, $\widetilde{\tau}\geq\tau$. This proves the existence results of Proposition \ref{strong3}. 

\bigskip

We now turn to the pathwise uniqueness of the solution. In the case   $\,\sigma_2^2-\sigma_1^2\geq\sigma_3^2-\sigma_2^2\,$, it follows from Theorem \ref{triple_coll_thm} that $\tau=\infty$ with probabilty one. In this case, pathwise uniqueness of the solution constructed above is a consequence of Theorem 3.2 in \cite{Ch}. The latter states that, for a system of stochastic differential equations with time-independent coefficients, strong existence in the presence of weak uniqueness implies pathwise uniqueness. 

\smallskip

In the case   $\,\sigma_2^2-\sigma_1^2<\sigma_3^2-\sigma_2^2\,$, we recall from Theorem \ref{triple_coll_thm} that we have $\mathbb{P}(\tau<\infty)>0$. In this case, we let $(X(t),W(t))$, $0\leq t\leq \tau$ and $\,(\widehat{X}(t),W(t))$, $0\leq t\leq \widehat{\tau}\,$ be two strong solutions of the equation \eqref{sde3}, where $\,\widehat{\tau}\,$ is the first time of a triple collision in the particle system $\widehat{X}=(\widehat{X}_1,\widehat{X}_2,\widehat{X}_3)$. By enlarging the probability space if necessary, we can extend $X$ and $\widehat{X}$ to weak solutions of the equation \eqref{sde3}, defined on the whole time interval $[0,\infty)$. Then, an application of Theorem 3.1 in \cite{Ch} shows that the joint laws of the triples $(\tau,X([0,\tau]),W([0,\tau]))$ and $(\widehat{\tau},\widehat{X}([0,\widehat{\tau}\,]),W([0,\widehat{\tau}\,]))$ are the same. We can now proceed as in the proof of Theorem 3.2 in \cite{Ch} to deduce $\,\tau=\widehat{\tau}\,$ and $X([0,\tau])=\widehat{X}([0,\widehat{\tau}\,])$ with probability one. \ep 

\subsection{Systems with finitely many particles}

We now turn to the proof of Theorem \ref{strong_sol} with $I=\{1,\dots,n\}$, where $n>3$. Although the main idea behind the construction of the strong solution is the same as for $n=3$, the proof is more involved here due to a more complicated pattern of collisions. For example, even in the absence of triple collisions, it is still   possible to have   collisions of the form $\,X_1(t)=X_2(t)\,$, $\,X_3(t)=X_4(t)\,$. 

\bigskip

\noindent\textbf{Proof of Theorem \ref{strong_sol} for $I=\{1,\dots,n\}$:} \textbf{Step 1.} As in the proof of Proposition \ref{strong3}, we start by constructing a strong solution to the system \eqref{sdeI} in an inductive manner. However, this time several layers of inductive constructions will be necessary. 

First, we recall the notation $Z^{[a,b],B,W}(b_1,b_2,c_1,c_2)$ for the strong solution of the system \eqref{sdeI} with $I=\{1,2\}$ and parameters $\,b_1,b_2,c_1,c_2\,$, which is driven by the independent standard Brownian motions $B$, $W$ on the time interval $[a,b]$. The latter exists by the results of section 5 in \cite{IFKP}. Next, we define a sequence of stopping times $0=\tau_0\leq\tau_1\leq\ldots\,$, subsets $A_0,A_1,A_2,\ldots$ of the set $\nn$ with $A_0=\emptyset$, and the desired strong solution of \eqref{sdeI} on $[\tau_0,\tau_1],[\tau_1,\tau_2],\ldots$ inductively by
\begin{eqnarray*}
X_{\pi_k(j)}(t)=X_{\pi_k(j)}(\tau_k)+\delta_{\pi_k(j)}(t-\tau_k)+\sigma_{\pi_k(j)} (W_{\pi_k(j)}(t) - W_{\pi_k(j)}(\tau_k)),\\
j:\,\{j-1,j\}\cap A_k=\emptyset,\;t\in[\tau_k,\tau_{k+1}],\\
\quad\\
X_{\pi_k(j)}[\tau_k,\tau_{k+1}]=\big(Z^{[\tau_k,\tau_{k+1}],W_{\pi_k(j)},W_{\pi_k(j+1)}}(\delta_j,\delta_{j+1},\sigma_j,\sigma_{j+1})\big)_1,\\
X_{\pi_k(j+1)}[\tau_k,\tau_{k+1}]=\big(Z^{[\tau_k,\tau_{k+1}],W_{\pi_k(j)},W_{\pi_k(j+1)}}(\delta_j,\delta_{j+1},\sigma_j,\sigma_{j+1})\big)_2,\\
j\in A_k,
\end{eqnarray*}
\begin{eqnarray*} 
&&\tau_{k+1}=\inf\{t>\tau_k:\;X_{(j)}(t)=X_{(j+1)}(t)\;\text{for some}\;j\notin A_k\},\\
&&A_{k+1}=\big\{1\leq j\leq n-1:\;X_{(j)}(\tau_{k+1})=X_{(j+1)}(\tau_{k+1})\big\}\,,
\end{eqnarray*}
where for each $k\geq0$, $\pi_k$ is a permutation of the set $\{1,\dots,n\}$ for which
\[
X_{\pi_k(1)}(\tau_k)\leq X_{\pi_k(2)}(\tau_k)\leq\ldots\leq X_{\pi_k(n)}(\tau_k).
\]
As in the case of $n=3$ it is not hard to see that this defines a strong solution of the system \eqref{sdeI} up to the time
\[
\tau^{\,[1]}\,:=\,\lim_{k\rightarrow\infty}\tau_k\,.
\]
We note at this point that, for $n>3$, $\tau^{[1]}=\tau$ is not necessarily true.  

\bigskip
\noindent\textbf{Step 2.} To proceed, for $1\leq j\leq n-1$ we let $(\rho_{j,k})_{k\geq1}$ be a (possibly empty) subsequence of the sequence $(\tau_k)_{k\geq1}$, which contains all the elements of the sequence $(\tau_k)_{k\geq1}$ for which $X_{(j)}(\tau_k)=X_{(j+1)}(\tau_k)$. Since for every $k_0\geq0$ there are at least two sequences of the form $(\rho_{j,k})_{k\geq1}$, which contain at least one element of $\{\tau_{k_0},\tau_{k_0+1}\}$, at least two of the sequences $(\rho_{j,k})_{k\geq1}$ are infinite. Thus, due to the continuity of the paths of the ordered particles, there exist $1\leq j_1<j_2\leq n-1$ such that
\[
X_{(j_1)}(\tau^{[1]})=X_{(j_1+1)}(\tau^{[1]}),\quad X_{(j_2)}(\tau^{[1]})=X_{(j_2+1)}(\tau^{[1]})
\]  
holds.

\medskip
If $j_2=j_1+1$, we have a triple collision of $X_{(j_1)}$, $X_{(j_1+1)}$ and $X_{(j_1+2)}$ at time $\tau^{[1]}$. In this case, $\tau\leq\tau^{[1]}$, and the existence results of the theorem readily follow.  

\medskip
If $\,j_2\neq j_1+1\,$, then we proceed with the construction of the processes $X_1,\ldots,X_n$ as in Step 1, but now starting at $t=\tau^{[1]}$ instead of $t=0$. This gives us a new sequence of stopping times $\tau^{[1]}=\tau^{[1]}_0\leq\tau^{[1]}_1\leq\ldots$ defined analogously to the stopping times $\tau_0\leq\tau_1\leq\ldots$ in Step 1. Next, we set $\tau^{[2]}=\lim_{k\rightarrow\infty}\tau^{[1]}_k$ and observe that we have constructed a strong solution to the system \eqref{sdeI} up to the time $\tau^{[2]}$. 

\medskip

Again, either there is a triple collision at $\tau^{[2]}$ and the proof is complete, or we proceed with the construction to extend the solution up to a time $\tau^{[3]}$. Proceeding in the same manner and assuming that we do not encounter a triple collision, we end up with a strong solution up to time $\tau^{[\infty]}:=\lim_{k\rightarrow\infty}\tau^{[k]}$. 

\medskip

Now, assuming that a triple collision has not occured, we conclude that, for each $k\in\nn$, there exist $1\leq j_1(k)<j_2(k)\leq n-1$ such that $j_2(k)>j_1(k)+1$ and 
\[
X_{(j_1(k))}(\tau^{[k]})=X_{(j_1(k)+1)}(\tau^{[k]}),\quad X_{(j_2(k))}(\tau^{[k]})=X_{(j_2(k)+1)}(\tau^{[k]}).
\] 
Arguing as before, we conclude that there exist $1\leq j_1(\infty)<j_2(\infty)<j_3(\infty)<j_4(\infty)\leq n-1$ such that
\[
X_{(j_\ell(\infty))}(\tau^{[\infty]})=X_{(j_\ell(\infty)+1)}(\tau^{[\infty]}),\quad \ell=1,2,3,4.
\]

Again, either there is a triple collision at time $\tau^{[\infty]}$, or we proceed with the construction of the strong solution until a time $\tilde{\tau}$ at which
\[
X_{(\,\widetilde{j}_\ell)}(\tilde{\tau})=X_{(\,\widetilde{j}_\ell+1)}(\tilde{\tau}),\quad \ell=1,\ldots 2^K,
\]
where $K$ is such that $\,2^K>n/2\,$ (unless a triple collision occurs before the time $\tilde{\tau}$). Now, the inequality $\,2^K> n/2\,$ implies that at time $\tilde{\tau}$ there must be a triple collision, which in turn implies $\tau\leq\tilde{\tau}$. Thus, we have constructed a strong solution up to the time $\tau$, as desired. 

\smallskip

Finally, pathwise uniqueness up to time $\tau$ can be shown by the same arguments as in the case $n=3$ (see the proof of Proposition \ref{strong3} for the details). \ep

\subsection{Systems with infinitely many particles}

We can now use the results of the previous subsection to construct the unique strong solution for the system with infinitely many particles.

\bigskip

\noindent{\bf Proof of Theorem \ref{strong_sol} for $\,I=\nn\,$:} Consider a probability space on which a system $W=(W_i:\;i\in I)$ of independent standard Brownian motions is defined. Then the result of the previous subsection shows that for every $n\in\nn$, there is a strong solution $X^{(n)}$ to the system \eqref{sdeI} with $I=\{1,\dots,n\}$ and the parameters $\delta_1,\dots,\delta_n,\sigma_1,\dots,\sigma_n$ which is defined up to the time
\eq
\inf\{t\geq0:\;X^{(n)}_i(t)=X^{(n)}_j(t)=X^{(n)}_k(t)\;\text{for\;some}\;1\leq i<j<k\leq n\}.
\en
In particular, it is defined for all $t\geq0$ if Condition \ref{cond2} is satisfied. Now, we can construct a strong solution $X$ of the system \eqref{sdeI} up to time $\tau$ by following the proof of Proposition 3.1 in \cite{Sh}, and making the following two modifications. We use the just described strong solutions $X^{(n)}$ instead of the weak solutions used in \cite{Sh} and replace the a priori estimate for the expected number of particles on intervals of the form $(-\infty,x]$ at a time $t\geq0$ in the proof of Proposition 3.1 in \cite{Sh} by the corresponding a priori estimate in the proof of our Proposition \ref{inf_syst}. This proves strong existence.

\bigskip 

To prove strong uniqueness, let $X'=(X'_i:\;i\in I)$ be another strong solution of the system \eqref{sdeI}, defined on the same probability space as $X$ and adapted to the same Brownian filtration, and let $\tau'$ be the corresponding first time of a triple collision. Moreover, we define the stopping times $\,0=\kappa'_0\leq\kappa'_1\leq\ldots\,$ and the sets $\,\{1,\dots,M\}=I'_0\subset I'_1\subset\ldots\,$ inductively by 
\begin{eqnarray*}
&&\kappa'_{k+1}=\Big(\inf\{t\geq\kappa'_k\, |\;\exists ~ i\in I'_k, j\notin I'_k:\;X'_j(t)=X'_i(t)\}\Big)\wedge\tau',\\
&&I'_{k+1}= \big\{ \, i\in I \,|\exists ~ 0\leq t\leq\kappa'_{k+1}:\;X'_i(t)=X'_j(t)\;\text{for\;some\;}j\in I'_k\, \big\}
\end{eqnarray*}
and let $0=\kappa_0\leq\kappa_1\leq\dots$ and $\{1,\dots,M\}=I_0\subset I_1\subset\dots$ be the corresponding quantities for the strong solution $X$. Then, by the strong uniqueness for the finite system, for each $k\in\nn$ the process
\[
\big(X'_1(t\wedge\kappa'_k),\dots,X'_{|I'_k|}(t\wedge\kappa'_k)\big)\,,\;\;t\geq0
\] 
must be the unique strong solution of \eqref{sdeI} with $I=I'_k$ and parameters $\delta_1,\dots,\delta_{|I'_k|}$, $\sigma_1,\dots,\sigma_{|I'_k|}$ driven by $W_1,\dots,W_{|I'_k|}$, stopped at $\kappa'_k$. Moreover, we have
\[
X'_i(t\wedge\kappa'_k)=X'_i(0)+\delta_i (t\wedge\kappa'_k) + \sigma_i W_i(t\wedge\kappa'_k)\,,\;\;t\geq0 
\] 
for all $i\notin I'_k$. The same arguments with $X'$ replaced by $X$ and induction over $k$ give
\eq
\kappa'_k=\kappa_k\,,\; I'_k=I_k\,,\; X(t\wedge\kappa_k)=X'(t\wedge\kappa'_k)
\en
for all $t\geq0$ and $k\in\nn$. Thus, by taking the limit $k\rightarrow\infty$ we obtain $\tau=\tau'$ and $X(t\wedge\tau)=X'(t\wedge\tau')$, $t\geq0$. 

Note that the almost sure identities $$\tau=\lim_{k\rightarrow\infty}\kappa_k \;\; \;\;\mathrm{and} \;\;\;\;\tau'=\lim_{k\rightarrow\infty}\kappa'_k$$ follow from the fact that, with probability one, for every $x\in\rr$ and $t\geq0$ there are finitely many particles on the interval $(-\infty,x]$ at time $t$ in $X$ and $X'$ (due to the a priori estimate in the proof of Proposition \ref{inf_syst} and the Borel-Cantelli Lemma). 
\ep

 \bigskip

\bibliographystyle{alpha}

\begin{thebibliography}{50}

\bibitem{BFK}
\textsc{Banner, A.}, \textsc{Fernholz, E.R.} \& \textsc{Karatzas, I.} (2005)
Atlas models of equity markets. \textit{Ann. Appl. Probab.} {\bf 15}, 2296-2330. 

\bibitem{BP}
\textsc{Bass, R.} \& \textsc{Pardoux, E.} (1987)
Uniqueness for diffusions with piecewise constant coefficients. \textit{Probab. Theory Related Fields} {\bf 76}, 557-572.

\bibitem{CL1}
\textsc{C\'epa, E.} \& \textsc{L\'epingle, D.} (1997)
Diffusing particles with electrostatic repulsion. \textit{Probab. Theory Related Fields} {\bf 107}, 429-449.

\bibitem{CL2}
\textsc{C\'epa, E.} \& \textsc{L\'epingle, D.} (2007)
No triple collisions for mutually repelling Brownian particles. \textit{Lecture Notes in Mathematics} {\bf 1899}, 241-246.

\bibitem{Ch}
\textsc{Cherny, A.S.} (2001)
On the uniqueness in law and the pathwise uniqueness for stochastic differential equations. \textit{Theory Probab. Appl.} \textbf{46}, 483-497. 

\bibitem{deB2}
\textsc{De$\,$Blassie, R.D.} (1998) On hitting single points by a multidimensional  diffusion. \textit{Stochastics} \textbf{65}, 1-11. 

\bibitem{deB}
\textsc{De$\,$Blassie, R.D.} (1999) Scale-invariant diffusions: transience and non-polar points. \textit{Bernoulli} \textbf{5}, 589-614. 

\bibitem{IFKP}
\textsc{Fernholz, E.R.}, \textsc{Ichiba, T.},  \textsc{Karatzas, I.} \& \textsc{Prokaj, V.} (2011)
Planar diffusion with rank-based characteristics, and a perturbed Tanaka equation. Available at \textit{http://arxiv.org/abs/1108.3992}.

\bibitem{Fr}
\textsc{Friedman, A.} (1974)
Nonattainability of a set by a diffusion process. \textit{Trans. Amer. Math. Soc.} {\bf 197}, 245-271.

\bibitem{Fr2}
\textsc{Friedman, A.} (2006)
\textit{Stochastic Differential Equations and Applications.} (Reprint of the Two-Volume    work published by Academic Press, 1975/76.) Two Volumes Bound as One, Dover Publications, Mineola, NY.

\bibitem{IK} 
\textsc{Ichiba, T.} \& \textsc{Karatzas, I.} (2010)
On collisions of Brownian particles. \textit{Ann. Appl. Probab.} \textbf{20}, 951-977.

\bibitem{IPBFK}
\textsc{Ichiba, T.}, \textsc{Papathanakos, V.}, \textsc{Banner, A.}, \textsc{Karatzas, I.} \& \textsc{Fernholz, E.R.}    (2011) Hybrid Atlas models. \textit{Ann. Appl. Probab.} {\bf 21},  609-644. 

\bibitem{Ram}
\textsc{Ramasubramanian, S.} (1983)
Recurrence of projections of diffusions. \textit{Sankhy$\bar{a}$ Ser. A} {\bf  45}, 20-31.

\bibitem{Ram2}
\textsc{Ramasubramanian, S.} (1988)
Hitting of submanifolds by diffusions. \textit{Probab. Theory Related Fields} {\bf  78}, 149-163.

\bibitem{RW}
\textsc{Reiman, M.} \& \textsc{Williams, R.J.} (1988)
A boundary property of semimartingale reflecting Brownian motions. \textit{Probab. Theory Related Fields} {\bf 77}, 87-97. 

\bibitem{Sh}
\textsc{Shkolnikov, M.} (2010)
Competing particle systems evolving by interacting L\'evy processes. \textit{Ann. Appl. Probab.}, to appear.  Available at \textit{http://arxiv.org/abs/1002.2811}.

\bibitem{VW}
\textsc{Varadhan S.R.S.} \& \textsc{Williams, R.J.} (1985)
Brownian motion in a wedge with oblique reflection. \textit{Comm. Pure Appl. Math.} {\bf 38} 405-443.

\bibitem{W1}
\textsc{Williams, R.J.} (1995)
Semimartingale reflecting Brownian motions in the orthant. In {\it ``Stochastic Networks''} (F.P. Kelly and R.J. Williams, eds.), Springer-Verlag.

\bibitem{W2}
\textsc{Williams, R.J.} (1987)
Reflected Brownian motion with skew symmetric data in a polyhedral domain. \textit{Probab. Theory Related Fields} \textbf{75}, 459-485.


\end{thebibliography}

\end{document}